\documentclass{conm-p-l}

\usepackage{eucal}

\begin{document}

\newcommand{\C}{{\mathbb C}}
\newcommand{\M}{{\mathcal M}}
\newcommand{\X}{{\mathcal X}}
\newcommand{\Oh}{{\mathbb O}}
\newcommand{\I}{{\mathbb I}}
\newcommand{\Gl}{{\rm Gl}}
\newcommand{\Hom}{{\rm Hom}}
\newcommand{\Fns}{{\rm Fns}}
\newcommand{\Cat}{{\mathcal C}}
\newcommand{\Q}{{\mathbb Q}}
\newcommand{\Z}{{\mathbb Z}}

\title{HKR characters and higher twisted sectors} \bigskip

\author{Jack Morava}
\address{Department of Mathematics, Johns Hopkins University, Baltimore,
Maryland 21218}
\email{jack@math.jhu.edu}
\thanks{The author was supported in part by the NSF}
\date {15 May 2004}

\begin{abstract} This is an introduction to work of Hopkins, Kuhn, and 
Ravenel [9] on generalized group characters, which seems to fit very well 
with the theory of what physicists call higher twisted sectors in the 
theory of orbifolds. \end{abstract}

\maketitle

\section{Basic definitions}

\noindent
This paper is meant to be expository. Its first main point is that 
the {\bf inertia stack} construction, which has been a focus of
considerable attention in recent work on orbifolds and stacks, can
be iterated to define a simplicial object $\I^\bullet(\X)$, which
is a convenient device for organizing the `higher twisted sectors' 
of the cohomology of an orbifold $\X$. The twisted sectors defined by 
the twice-iterated intertia stack construction are crucially important 
in elliptic cohomology ([8], cf. also [16]), but the significance in 
physics of even higher iterations is not yet clear. They are nevertheless 
interesting invariants, and my second main point is that these things 
{\bf already} have a deep literature in mathematics. The study of conjugacy
classes of groups is absolutely fundamental to the theory of linear
representations, and it seems that conjugacy classes of commuting {\bf
tuples} of group elements play a systematically analogous role in the
theory of actions on manifolds. \bigskip

\noindent
This section recalls some standard facts from the `classical' theory
of orbifolds. The second section is concerned with the semisimplicial
constructions mentioned above, and the third is a quick account of the
higher character theory of Hopkins, Kuhn, and Ravenel. \bigskip

\noindent
This paper began as a talk at the ChengDu (Sichuan) satellite ICM conference 
on stringy orbifolds. I want to thank Alejandro \'Adem and Yongbin Ruan on
one hand, and Arkady Vaintrob and the referee on the other, for their
interest and helpful comments. I owe Ian Leary a special acknowledgement, 
for contributing an important idea at a crucial point (\S 2.3). This paper
grew out of many conversations with Matthew Ando; in a better world he would
be its coauthor. \bigskip

\noindent
{\bf 1.0} I will work in a convenient {\it ad hoc} category of orbispaces. 
For my purposes, an orbispace is a (topological) category $\X := [X/G]$ 
(with the points of a topological space $X$ as objects), defined by an 
action of a compact Lie group $G$ on $X$, subject to the restriction that 
the isotropy group $G_x$ of any point $x \in X$ be {\bf finite}. Morphisms 
of orbispaces are to be {\bf equivalence classes}, under invertible natural 
transformations, of (continuous) functors between categories. This class 
is rich enough to contain some interesting examples: \medskip

\noindent
{\bf Ex 1} If $\M$ is a reduced $d$-dimensional orbifold, then its principal 
orthogonal frame `bundle' $\Oh(M)$ is a smooth manifold upon which the 
orthogonal group $\Oh(d)$ acts with finite isotropy. By a fundamental lemma 
[10 \S 1 (example)] of Kawasaki, the category (or groupoid) $[\Oh(M)/\Oh(d)]$ 
is equivalent to the category defined by the original orbifold $\M$.\medskip

\noindent
{\bf Ex 2} If $G$ is a {\bf finite} group, then the category $[*/G]$ with one
object, and the set $G$ of morphisms, is an interesting {\bf un}reduced 
orbifold. \bigskip

\noindent
{\bf Remarks:} Useful topological constructions take us out of the category 
of smooth objects, so it is convenient to work with a class slightly larger 
than the usual orbifolds. In general, I will use the {\tt mathcal} typeface 
for an orbispace, and the usual mathematical typeface for its underlying 
space of objects; thus $\X := [X/G]$ has objects $X$ and underlying quotient 
space $X/G$. However, there will be exceptions: \bigskip 

\noindent
{\bf 1.1} If $G$ is a group, and $X \in (G-{\rm spaces})$, then
\[
I(X) := \{(g,x) \in G \times X \;|\: gx = x \}
\]
is itself a $G$-space, with action defined by 
\[
h \cdot (g,x) = (hgh^{-1},hx) \;.
\]
$I$ is thus a functor from the category of $G$-spaces to itself. The isotropy
group of $(g,x) \in I(X)$ is 
\[
\{ h \in G \;|\; h(g,x) = (hgh^{-1},hx) = (g,x) \} \;;
\]
being a subset of $G_x$, it is finite if the latter is. It follows that if 
$\X = [X/G]$ is an orbispace, in the sense above, then
\[
\I(\X) := [I(X)/G] 
\]
is again such an orbispace; following the terminology of algebraic
geometers, it is now called the {\bf inertia stack} of $\X$. [It is
also the fixed-point orbispace [10] of the circle group, acting on the
free loops in $\X$.] The description above makes it clear 
that $\I$ is an endofunctor of the category of orbispaces. \bigskip

\noindent
{\bf 1.2} These constructions define some useful invariants. I will call 
the Borel cohomology 
\[
H^*(\X,\Q) := H_G^*(X,\Q) := H^*(EG \times_G X,\Q)
\]
(in this paper all coefficients will be vector spaces over $\Q$) the 
{\it ordinary} cohomology of the orbispace $\X$: its Leray spectral
sequence has as $E_2$-term, the cohomology
\[
H^*(X/G,{\mathcal H}^*(G_x,\Q)) 
\]
of the quotient with coefficients in a sheaf whose stalk at $x$ is the group 
cohomology of $G_x$. Since these groups are by hypothesis finite, this 
sheaf is concentrated in degree zero, and the spectral sequence degenerates 
to an isomorphism
\[
H^*(\X,\Q) \cong H^*(X/G,\Q) = H^0(G,H^*(X,\Q))
\]
with the cohomology of the quotient space. This is interesting enough, but 
it is not very subtle. \bigskip

\noindent
A more powerful invariant is defined by the equivariant $K$-theory
\[
K^*(\X) := K_G^*(X)
\]
of the orbispace. \bigskip

\noindent
{\bf 1.3 Theorem:} There is a natural multiplicative transformation
\[
K^*(\X) = K_G(X)^* \to H^*_G(I(X),\Q) = H^*(\I(\X),\Q)
\]
which becomes an isomorphism after tensoring with $\Q$ on the left. 
\bigskip

\noindent
{\bf Remarks:} This is proved in [2]; nowadays this (rational) 
invariant is usually called the Adem-Ruan, or classical, orbifold cohomology.
It is to be distinguished from the Chen-Ruan [4] orbifold cohomology, which 
has a different multiplicative structure. [I will ignore some deep questions 
about the gradings of these theories, since I have nothing to say about them.]
I should note that the existence of some such generalized Chern character 
was also known to Baum and Brylinski [3]. \bigskip

\noindent
The cohomology groups on the right have a natural decomposition [2 \S 5.1] as
\[
\oplus_{g \in \hat G} \; H^*_{C(g)}(X^g,\Q)
\]
where $\hat G$ denotes the set of conjugacy classes in $G$, $X^g$ is the
set of $g$-fixed points in $X$, and $C(g)$ is the centralizer of $g$ in
$G$. [More precisely: for any choice of $g$ in the appropriate conjugacy
class, the cohomology group in question is well-defined under conjugation 
by elements of $G$.] \bigskip

\noindent
The contributions to this sum, indexed by conjugacy classes {\bf other
than the identity} are now called the {\bf twisted sectors} of the 
cohomology. \bigskip

\section{Higher inertia stacks} 

\noindent
{\bf 2.1 Definition:} If $\X = [X/G]$ as above, let
\[
\I^n(\X) = [I^n(X)/G] \;;
\]
note that 
\[
I^n(X) = \{ (g_1,\dots,g_n;x) \in G^n \times X \;|\; g_i \in G_x, \; \forall
i,k \; [g_i,g_k] = 1 \} \;.
\]
{\bf Proof:} See the argument in \S 1.1, and induct. \bigskip

\noindent
{\bf For example:} If $X = *$ is a single point, 
\[
\I^n[*/G] = \Hom(\Z^n,G)/G
\]
is the set of conjugacy classes of commuting $n$-tuples of elements in
$G$. When $n=1$, this is just the classical set $\hat G$ of conjugacy classes 
in $G$. \bigskip

\noindent
The construction of the inertia stack is essentially local, so
more generally
\[
\I^n[X/G] = [(\bigcup_{x \in X} I^n[\{x\}/G_x])/G] \;.
\] \bigskip

\noindent
{\bf 2.2} Recall now that a simplicial object in a category $\Cat$
can be defined as a functor $C$ from the category of finite ordered sets to
$\Cat$. We can think of such a functor as defined by its sets $C[n]$
of $n$-simplices, together with various face and degeneracy maps between
them. \bigskip

\noindent
A simplicial object in the category of spaces (for example, a simplicial set) 
has a {\bf geometric realization}
\[
|C| = \coprod_{n \geq 0} (C[n] \times \Delta^n)/({\rm face \; \& \; degeneracy
\; relations}) \;.
\]
For example: a category $\Cat$ can be regarded, following Grothendieck and
Segal [17], as a simplicial set with objects as zero-simplices, morphisms as 
one-simplices, and chains of $n$ composable morphisms as its $n$-simplices. 
The face maps are defined by composing maps, and degeneracies are defined
by inserting identities. The geometric realization of this simplicial set 
is sometimes called the classifying space for the category; in particular, 
\[
|[*/G]| = BG 
\]
is the classifying space for the (discrete) group $G$, and more generally
the geometric realization 
\[
|[X/G]| = EG \times_G X
\]
of (the category defined by) a transformation group is homotopy equivalent to 
its Borel construction. The map 
\[
|[X/G]| = EG \times_G X \to * \times_G X = X/G
\]
which collapses (the free contractible $G$-space) $EG$ to a point is
sometimes called the `homotopy-to-geometric' quotient; the arguments
of \S 1.2 above show that for our class of orbispaces, this map induces
an isomorphism on rational cohomology. \bigskip

\noindent
{\bf 2.3} Conjugation by a group element defines a functor from the category 
$[*/G]$ to itself; thus $G$ acts on $|[*/G]| = BG$. This endofunctor is 
naturally equivalent to the identity: natural transformations of functors 
become (un-base pointed) homotopies under geometric realization, but the 
resulting action of $G$ on $BG$ need not be trivial. This construction 
generalizes to an action of $G$ on the category $[X/G]$ and its realization 
$|[X/G]|$: thus the homotopy quotient of a $G$-action is again a (not 
necessarily trivial!) $G$-space. I will write $X//G$ for the quotient of the 
homotopy quotient by $G$: the homotopy-to-geometric quotient map can thus be 
factored, as
\[
|[X/G]| \to |[X/G]|/G := X//G \to X/G \;.
\]
I owe thanks to Ian Leary for help in understanding this intermediate
quotient. He notes that if $G \times G$ acts on $G$ by $(g,h) \cdot k \mapsto
gkh^{-1}$ then the usual construction for $EG$ (as a simplicial set with
$n$-simplices $G^{n+1}$) becomes a space with $G \times G$-action; as 
such it is a classifying space for the family consisting of subgroups
of $G \times G$ which are conjugate to a subgroup of the diagonal. It 
follows that
\[
EG/(G \times G) = BG/G = *//G \;,
\]
and that $\pi_1(*//G) \cong H_1(G,\Z)$ is the abelianization of $G$;
the quotient $|[*/G]| \to *//G$ is a homotopy equivalence if and only
if $G$ is abelian. \bigskip

\noindent
{\bf 2.4} The simplicial set $n \mapsto \Z^n$ defining $[*/\Z]$ is in fact
a simplicial object in the category of abelian groups: group composition
is a homomorphism when the group is abelian. It follows that the
{\bf co}variant functor
\[
n \mapsto (\check{\Z})^n := \Hom(\Z^n,\Z)
\]
is, in a natural sense, a {\bf co}simplicial abelian group. 
\bigskip

\noindent
{\bf Definition:} The functor 
\[
n \mapsto \Hom((\check{\Z})^n,G)/G
\]
defines the simplicial set $\I^\bullet[*/G]$ of conjugacy classes of 
commuting tuples of elements $G$, cf. [11 \S 4]; more generally,
\[
\I^\bullet[X/G] := [(\bigcup_{x \in X} \I^\bullet[\{x\}/G_x])/G]
\]
is the {\bf simplicial inertia stack} of $\X$. \bigskip

\noindent
We can use this construction to elaborate Adem and Ruan's construction for
orbifold cohomology:
\[
n \mapsto H^*(\I^n(\X),\Q)
\]
is a cosimplicial object in the category of graded-commutative algebras,
which keeps simultaneous track of the higher inertia stacks of $\X$. 
\bigskip

\noindent
{\bf 2.5 Theorem:} There is a natural transformation
\[
|\I^\bullet[X/G]| \to X//G 
\]
which is an equivalence if $G$ is abelian. \bigskip

\noindent
The proof is by construction; it is easiest to begin in the special
case when $X$ is a point. Then $\I^\bullet[*/G]$ is a simplicial set 
with one zero-simplex; a one-simplex is a conjugacy class, a two-simplex
is a conjugacy class of commuting pairs, etc. If $\langle g_1,\dots,g_n 
\rangle$ is an $n$-simplex, then its faces are the maps
\[
\langle g_1,\dots,g_n \rangle \mapsto \langle g_1,\dots,g_{i-1}g_i,\dots, g_n 
\rangle
\]
and its degeneracies are the maps which insert identity elements. These
are exactly the maps defining the classifying space of $G$; but we are
working now not with group elements, but conjugacy classes. \bigskip

\noindent
The promised map is then the quotient of the obvious equivariant inclusion
\[
\Hom((\check{\Z})^\bullet,G) \to BG   
\]
by $G$. Because $\I^\bullet$ is a local construction, this definition now
extends directly to $[X/G]$; alternately, we can display the simplicial 
object $\I^\bullet[X/G]$ (with most of its maps supressed) as
\[
\dots \to \coprod [(X^g \cap X^h)/C(g,h)] \to \coprod [X^g/C(g)] \to [X/G] \;,
\]
where the $n$th coproduct is indexed by conjugacy classes of commuting
$n$-tuples, and $C(g_1,\dots,g_n)$ is the centralizer of the commuting
tuple. \bigskip

\noindent
{\bf Remarks:} It is tempting to think of this construction as a kind of 
blowup or resolution of the Borel construction; it seems analogous in 
some ways to Segal's [17] reconstruction of a manifold, up 
to homotopy, from the category defined by the sets of an atlas with
inclusions as morphisms. In our case, the charts are reminiscent of the 
complete sets of commuting observable of classical quantum mechanics. 
Kuhn [11 \S 7] remarks that $|\I^\bullet[*/G]|$ is in fact a $\Gamma$-space 
[though not, in general, a special $\Gamma$-space] in the sense of 
Segal.\bigskip

\noindent
I am reluctant to admit that I don't know how a single example works out. 
Symmetric groups and finite subgroups of ${\rm Sl}_2(\C)$ are of course
very interesting candidates. \bigskip

\noindent
This construction may also be related to the theory of motivic
integration: if $X$ is an algebraic variety, say over the complexes,
the $n$-simplices of $\I^\bullet$ are roughly deformations of the scheme 
over fields of transcendence degree $n$. To make this precise 
would require a better understanding of the degree-shifting numbers 
[4; 12 \S 8; 14 \S 2], which do not appear in the formalism above. \bigskip

\noindent
When $n=1$, these are locally constant $\Q$-valued functions $w$ on the 
fixed-point set $X^g$, which are slightly more sophisticated than the
function which assigns to $g$, the number
\[
\log \; \det \; (g|\nu) \;,
\]
where $g|\nu$ represents the action of $g$ on the normal bundle of 
$X^g$ in $X$. [I'm assuming here that the orbifolds in question 
have complex structures on their tangent spaces.] In general, the 
normal bundle to the fixed point set of a commuting $n$-tuple 
$\langle g_1,\dots,g_n \rangle$ has a flag decomposition as a sum 
\[
X^{g_1} \cap \dots \cap X^{g_{i-1}} \subset X^{g_1} \cap \dots \cap X^{g_i} 
\]
of normal bundles, and it seems reasonable to expect that the degree-shifting 
number of this $n$-tuple will be the sum of the degree-shifting numbers of 
these subbundles. \bigskip

\section{HKR characters}

\noindent
{\bf 3.0} A homomorphism from a free abelian group to a finite group $G$ 
factors through some finite abelian quotient group, so
\[
\Hom(\Z^n,G)/G = \Hom(\hat{\Z}^n,G)/G = {\prod}^*_p \Hom(\Z_p^n,G)/G
\]
decomposes as the restricted product (with only finitely many nontrivial 
entries) of $p$-local contributions, indexed by primes $p$. This uses the
fact that
\[
\hat \Z = {\prod}_p \Z_p \;,
\]
where $\Z_p = \lim \Z/p^n \Z$ the $p$-adic integers. \bigskip

\noindent
Since products of simplicial sets (and spaces) are defined coordinate-wise, 
$\I^\bullet [X/G]$ can be expressed as a restricted infinite fiber
product (over $[X/G]$) of $p$-local objects $ \I^\bullet_p[X/G]$ built 
like $\I$ but with $\Z_p$ replacing $\Z$. I will ignore questions
about infinite restricted products by assuming that $[X/G]$ is `ramified'
at a finite set of primes (dividing $\# G$, say, when the group is finite);
the cohomology of the simplicial inertia stack can then be calculated 
from the local contributions, one prime at a time. \bigskip

\noindent
{\bf 3.1} In this context, Hopkins, Kuhn, and Ravenel [9] provide, for 
each $n \geq 1$, an interpretation of $H^*(\I^n_p[X/G],\Q)$ which reduces 
when $n = 1$ to the theorem of Adem and Ruan in \S 1.3 above. To state 
these results, however, requires a short digression about cobordism. 
\bigskip

\noindent
Very briefly, then: cobordism is to homology as smooth manifolds are
to simplices. In this theory, a $d$-dimensional chain is not some 
sum of nasty singular simplices, but a map, say $f: M \to X$, of a nice 
smooth $d$-manifold $M$ to the space $X$ of interest. Instead of boundaries 
of simplices, we take boundaries of manifolds; thus $\partial f : \partial M 
\to X$ is the boundary of $f$, which is said to be closed if $\partial M = 0$.
Similarly, $f = \partial F$ if $\exists F: W \to M$ such that $\partial
W = M$ and $F|_{\partial W} = f$. The analog of the homology of $X$ is
the quotient of the abelian semigroup of closed objects (cycles) by the
subsemigroup of boundaries; this is well-defined, since of course $\partial 
\circ \partial = 0$. It is more usual to say that these groups
are defined by classes of maps of smooth manifolds to $X$ under the 
equivalence relation defined by cobordism: which is to say that two maps
of closed manifolds to $X$ are related if they are the boundary values 
of maps defined on a smooth manifold of one higher dimension. \bigskip

\noindent
These groups are obviously homotopy-invariant (use the cobordism defined
by a cylinder) and covariant: a map $\phi : X \to Y$ pushes the class $[f]$
to the class $[\phi \circ f]$. Atiyah's convention is to call this 
(graded-abelian-group-valued, homological) functor the {\bf bordism} of $X$;
there is a corresponding {\bf co}homological theory (contravariant under
pullback or fiber product, using Thom's theory of transversality), now 
usually called the {\bf co}bordism of $X$. One advantage of the latter 
theory is a nice multiplicative structure, defined by the obvious Cartesian
product, without need for any Eilenberg-Zilber foolishness. \bigskip

\noindent
Cobordism theory has very natural connections with the theory of group
actions on manifolds: the Borel construction 
\[
EG \times_G M \to EG \times_G * = BG
\]
associated to a $G$-manifold $M$ is a kind of relative manifold, which
defines a $(-d)$-dimensional class in the cobordism of $BG$. This is
the beginnings of a rich subject; a more sophisticated approach can be 
found in [7]. HKR theory is a natural generalization of the classical 
theory of characters of representations of groups on vector spaces to 
a theory of characters for actions on manifolds. \bigskip

\noindent
The advantages of cobordism (geometric naturality, etc.) are recognized
in the Russian literature, where it is usually called `intrinsic homology'. 
Its {\bf dis}advantages include the fact that there are many cobordism 
theories, depending on one's favorite choice of manifold: oriented, spin, 
symplectic, framed \dots each with its own special features. A more 
substantial issue is that the ground ring of such a theory (ie, the value 
of the cohomology theory on a point) tends to be quite large. It is 
arguably the cobordism theory of stably almost complex manifolds (with 
a complex structure on the sum of the tangent bundle with some trivial 
bundle) which is technically most accessible; that theory, called complex 
cobordism, has a polynomial ground ring 
\[
MU^* := MU^*({\rm pt}) \cong MU^* \cong \Z[x_i \;|\; i \geq 1]
\]
with one generators of each even degree. [Frank Adams's convention is to 
write $ML^*$ for the cobordism theory of manifolds with structure group 
reduced to the Lie group $L$, eg $U$ for weakly almost complex manifolds]. 
Over the rationals,
\[
MU^*({\rm pt}) \otimes \Q = \Q [\C P_n \;|\;n \geq 1]
\]
is the polynomial ring generated by the complex projective spaces; but 
these classes do not generate over the integers. \bigskip

\noindent
More generally, an old argument of Dold shows that there is
a natural multiplicative transformation
\[
MU^*(X) \to H^*(X, MU^* \otimes \Q)
\]
which factors through an isomorphism of the rationalization of the 
left-hand side. Over the rationals, then, there is in some sense little
difference between cobordism and ordinary cohomology. [In some sense
the chromatic filtration and the Atiyah-Swan filtration of ordinary
cohomology by size of supporting elementary abelian subgroups are aspects
of some common underlying phenomenon; see [5,6] for recent work on
the latter topic.] The advantage of cobordism lies in its naturality: its 
cycles are geometric objects, which carry characteristic class data (for 
example, of the sort familiar to physicists in the theory of `gravitational 
descendents'). \bigskip

\noindent
{\bf 3.2} However, These cobordism theories are often too big to be 
technically convenient - for example, their ground rings are not 
Noetherian - so topologists have developed an arsenal of techniques to make 
them more useful. One useful ruse is to work $p$-locally, at some fixed prime.
It turns out that to understand $MU$ in general, it suffices to understand
a hierarchy of cohomology theories with ground rings
\[
\hat E^*_n = \Z_p [v_1,\dots,v_{n-1}]((v_n^{-1})) 
\]
indexed by integers $n \geq 1$, defined as truncations (in a suitable 
sense) of the $p$-completion of $MU$; here $v_k$ can be taken to be
the cobordism class of a degree $p$ hypersurface in $\C P(p^k)$, and 
$A((x))$ is the formal Laurent series extension of a ring $A$ which
allows only finitely many negative powers of $x$. When $n=1$, this 
theory is a version of $p$-adically completed complex $K$-theory. 
\bigskip

\noindent
The study of these theories tends to involve some quite subtle number
theory, and one of the main technical advances in [9] is the construction
of a certain faithfully flat ring extension $\hat E_n \subset \hat D_n$,
which is most naturally interpreted as a kind of generalized Galois 
extension, with Galois group $\Gl_n(\Z_p)$. \bigskip

\noindent
{\bf Theorem:} There is a natural multiplicative transformation
\[
\hat E^*_n(|[X/G]|) \to H^*(\I^n_p[X/G],\hat D_n \otimes \Q)^{\Gl_n(\hat 
\Z_p)-{\rm inv}}
\]
which factors through an isomorphism with the rationalization of the
group on the left. \bigskip

\noindent
The term on the right is the subring of invariants under the action of
the Galois group $\Gl_n(\Z_p)$, but that action requires some 
clarification. The point is that this group acts on the coefficient
ring $\hat D_n$, but it also acts on $\I^n_p$, through its construction
in terms of conjugacy classes of homomorphisms from $\Z_p^n$ to $G$. 
The relevant action on the right is the (diagonal) product of these two
natural actions. These groups can be expressed in terms of fixedpoint
sets by formulae generalizing [2 \S 5.1]. \bigskip

\noindent
{\bf 3.3} This indeed restricts when $n=1$ to the theorem of \S 1.3, plus
(a $p$-adic version of) a theorem of Artin: there is a natural multiplicative
transformation 
\[
R(G) \to \Fns(\hat G,\Q_{cyc})^{{\rm Gal}(\Q_{cyc}/\Q)-{\rm inv}}
\]
which factors through an isomorphism with the rationalization of the
left-hand side; where $\Q_{cyc}$ is the cyclotomic closure of the rationals
(defined by adjoining all roots of unity). \bigskip

\noindent
This natural transformation is nothing but the map which assigns to 
a representation, its character; this version of the theorem encompasses
the fact, also due to Artin, that the values of such characters lie
in $\Q_{cyc}$. The Galois group 
\[
{\rm Gal}(\Q_{cyc}/\Q) \cong \hat \Z^\times
\]
of this extension is the multiplicative group of profinite integers,
whose $p$-local component is the $p$-adic unit group
\[
\Z_p^\times = \Gl_1(\Z_p) \;.
\]
The statement above conceals an action of $\hat \Z^\times$ on the conjugacy
classes, in which $k \in \Z$ sends the class of $g$ to the class of $g^k$
(away from the order of $g$). \bigskip

\noindent
{\bf 3.4} Here are a few closing remarks: \medskip

\noindent
i) The rings $\hat E_n$ classify (in a suitable sense) one-dimensional formal
groups of height $n$ over $p$-adic integer rings, and the rings $\hat D_n$ 
classify such groups, together with a level structure: this is a preferred
basis for the torsion subgroup. \medskip

\noindent
In the theory of algebraic stacks [1], the cyclotomic Galois action 
plays a distinguished role; the level structure is just a choice of 
isomorphism of $\Q_p/\Z_p$ with the group of $p$-power roots of unity. In 
the case of a stack defined over $\Q$, it is natural to think of the center 
of $\Gl_n(\hat Z_p)$ as acting through the determinant 
\[
\det : \Gl_n(\Z_p) \to \Z_p^\times
\]
on the roots of unity. \medskip

\noindent
ii) The $\hat E_n$'s and the $\hat D_n$'s do {\bf not} fit together
naturally as a (co)simplicial ring. In particular, the natural 
action of the symmetric group $\Sigma_n$ on $\I^n$ gets lost in the
action of $\Gl_n$ on $\hat D_n$. \medskip

\noindent
This suggests that there is lots of room in the transition between
chromatic levels for all sorts of gerbish orbifold twisting, and other
kinds of noncommutative monkey business \dots \bigskip

\newpage

\bibliographystyle{amsplain}

\begin{thebibliography}{99}

\bibitem [1]{1} D. Abramovich, T. Graber, A. Vistoli, Algebraic orbifold
quantum products, in {\bf Orbifolds in mathematics and physics} 1 - 24,
Contemporary Math 310, AMS (2002); available at {\tt math.AG/00112004}

\bibitem [2]{2} A. Adem, Y.B. Ruan: Twisted orbifold $K$-theory, Comm. 
Math. Phys. 237 (2003) 535 - 556; available at {\tt math.AT/0107168}

\bibitem [3]{3} P. Baum, J.L. Brylinski: Noncommutative topology: talk at 
AMS Winter meeting (2000)

\bibitem [4]{4} W. Chen, Y.B. Ruan, Orbifold quantum cohomology, in {\bf 
Orbifolds in mathematics and physics} 25 - 85, Contemporary Math 310, 
AMS (2002); available at {\tt math.AG/0005198}

\bibitem [5]{5} D.J. Green, I. Leary, The spectrum of the Chern subring, 
Comment. Math. Helv. 73 (1998) 406 - 426

\bibitem [6]{6} ----------, J. Hunton, B. Schuster, Chromatic characteristic 
classes in ordinary group cohomology, Topology 42 (2003) 243 - 263

\bibitem [7]{7} J. P. C. Greenlees, N. P. Strickland: Varieties and local 
cohomology for chromatic group cohomology rings, Topology 38 (1999) 
1093 - 1139

\bibitem [8]{8} M. J. Hopkins: Characters and elliptic cohomology, in 
{\bf Advances in homotopy theory} (Cortona, 1988), 87 - 104; London Math. 
Soc. Lecture Notes 139, Cambridge Univ. Press (1989)

\bibitem  [9]{9} --------, N. J. Kuhn, D. C. Ravenel: Generalized group 
characters and complex oriented cohomology theories. J. Amer. Math. Soc. 13 
(2000) 553 - 594

\bibitem [10]{10} T. Kawasaki, The signature theorem for $V$-manifolds,
Topology 17 (1978) 75 - 83

\bibitem [11]{11} N. J. Kuhn: Character rings in algebraic topology, in 
{\bf Advances in homotopy theory} (Cortona, 1988), 111 - 126, London Math. 
Soc. Lecture Notes 139, Cambridge Univ. Press (1989)

\bibitem [12]{12} E. Looijenga: Motivic measures, in {\bf S\'eminaire 
Bourbaki}, Asterisque 276 (2002)

\bibitem [13]{13} E. Lupercio, B. Uribe: Loop groupoids, gerbes, and 
twisted sectors on orbifolds, in {\bf Orbifolds in mathematics and physics} 
163 - 184 - 24, Contemporary Math 310, AMS (2002); available at 
{\tt math.AT/0110207}

\bibitem [14]{14} J. Morava: Some Weil group representations motivated by
algebraic topology, in {\bf Elliptic curves and modular forms in algebraic
topology}, Springer Lecture Notes in Mathematics 1326 (1986)

\bibitem [15]{15} M. Reid: La correspondance de McKay, in {\bf Seminaire 
Bourbaki}, Asterisque 276 (2002)

\bibitem [16]{16} S. Norton: Appendix [on generalized moonshine] to 
G. Mason, Finite groups and modular functions, in Proc. Sympos. Pure Math., 
47: {\bf Arcata Conference} 181--210 (1987)

\bibitem [17]{17} G. Segal: Classifying spaces and spectral sequences,
Publ. Math. IHES 34 (1968)

\end{thebibliography}

\end{document}